\documentclass[12pt]{amsart}

\usepackage{amssymb}

\usepackage{enumerate}

\usepackage{graphicx}

\makeatletter
\@namedef{subjclassname@2010}{%
  \textup{2010} Mathematics Subject Classification}
\makeatother

\newtheorem{thm}{Theorem}

\newtheorem{lem}[thm]{Lemma}

\newtheorem{con}[thm]{Conjecture}

\theoremstyle{definition}

\begin{document}


\baselineskip=17pt



\title[On sequences of integers with small prime factors]{On sequences of integers with small prime factors}

\author[C.L.~Stewart]{C.L.~Stewart}
\address{Department of Pure Mathematics, University of Waterloo \\
Waterloo, Ontario, Canada N2L 3G1}
\email{cstewart@uwaterloo.ca}

\dedicatory{For Professor Henryk Iwaniec on the occasion of his seventy-fifth birthday}

\date{}

\begin{abstract}
We show that the difference between consecutive terms in sequences  of integers whose greatest prime factor grows slowly tends to infinity.
\end{abstract}



\maketitle

\section{Introduction}

Let $y$ be a real number with $y\geq 3$ and let $1=n_1<n_2<n_3< ...$ be the increasing sequence of positive integers composed of primes of size at most $y$. In 1908 Thue \cite{Thue} proved that
\begin{equation}\label{eq1}
\lim_{i \to\infty} n_{i+1} -n_{i}=\infty,
\end{equation}
see also P\'olya \cite{P} and Erd\H{o}s \cite{Erd}. Thue's result was ineffective. In particular his proof does not allow one to determine, for every positive integer $m$, an integer $i(m)$ such that $n_{i+1} -n_{i}$ exceeds $m$ whenever $i$ is larger than $i(m)$.
Cassels \cite{Cassels} showed how \eqref{eq1} can be made effective by means of estimates due to Gelfond \cite{G} for linear forms in two logarithms of algebraic numbers. In 1973 Tijdeman \cite{Tijd} proved, by appealing to work of Baker \cite{Baker} on estimates for linear forms in the logarithms of algebraic numbers, that there is a positive number $c$, which is effectively computable in terms of $y$, such that 
\begin{equation} \label{2}
n_{i+1} -n_{i} > n_{i}/(\log n_i)^c
\end{equation}
for $n_i \geq 3$. In addition Tijdeman showed that there are arbitrarily large integers $n_i$ for which \eqref{2} fails to hold when $c$ is less than $\pi(y)-1$; here $\pi(x)$ denotes the counting function for the primes up to $x$.

Now let $y=y(x)$ denote a non-decreasing function from the positive real numbers to the real numbers of size at least $3$. For any integer $n$ let $P(n)$ denote the greatest prime factor of $n$ with the convention that $P(0)=P(\pm 1)=1$. Let $(n_i)^\infty_{i=1}$ be the increasing sequence of positive integers $n_i$ for which
\begin{equation} \label{3}
P(n_{i}) \leq y(n_{i}).
\end{equation}
For any integer $k$ with $k\geq 2$ let $\log_k$ denote the $k$-th iterate of the function $x \rightarrow \max(1,\log x)$ for $x>0$. We shall prove that \eqref{eq1} holds provided that
\begin{equation} \label{4}
 y(n)= o(\frac{\log_2 n\log_3 n}{\log_4 n}).
\end{equation}
Furthermore if we assume the abc conjecture, see \S 2, then we can prove that \eqref{eq1} holds provided that
\begin{equation} \label{5}
 y(n)= o(\log n).
\end{equation}

For any real number $x\geq 2$ put
$$
\delta(x) = \exp(\frac{x\log_2 x}{\log x}).
$$
 We shall deduce \eqref{4} from the following result.

\begin{thm} \label{Theorem 1}
Let $y=y(x)$ be a non-decreasing function from the positive real numbers to the real numbers of size at least $3$. 
 Let $(n_1,n_2,...)$ be the increasing sequence of positive integers $n_i$ for which \eqref{3} holds. There is an effectively computable positive number $c$ such that for $i\geq 3,$
 \begin{equation} \label{6}
n_{i+1} -n_{i} > n_i/ (\log n_i)^{\delta(cy(n_{i+1}))}.
\end{equation}
Furthermore there is an effectively computable positive number $c_1$ such that for infinitely many positive integers $i$
 \begin{equation} \label{7}
n_{i+1} -n_{i} <n_i\exp(c_1y(n_i))/ (\log n_i)^{r-1},
\end{equation}
where $r=\pi(y(\sqrt{n_i})).$
\end{thm}

Observe that we obtain \eqref{eq1} from \eqref{6} when \eqref{4} holds on noting that in this case $n_{i+1} \leq 2n_i$ and
$$
(\log n)^{\delta(cy(n))} = o(n).
$$
In order to establish \eqref{6} we shall appeal to an estimate for linear forms in the logarithms of rational numbers due to Matveev \cite{Matveev}, \cite{Matveev2}. The upper bound \eqref{7} follows from an averaging argument based on a result of Ennola \cite{Ennola}.

We are able to refine the lower bound \eqref{6} provided that the abc conjecture is true.

\begin{thm} \label{Theorem 2}
Let $y=y(x)$ be a non-decreasing function from the positive real numbers to the real numbers of size at least $3$. Let $(n_1,n_2,...)$ be the increasing sequence of positive integers $n_i$ for which \eqref{3} holds and let $\varepsilon$ be a positive real number.
If the abc conjecture is true then there exists a positive number $c_1=c_1(\varepsilon)$, which depends on $\varepsilon$, and a positive number $c_2$ such that for $i\geq 1$,
\begin{equation} \label{8}
n_{i+1} -n_{i} > c_1(\varepsilon)n_i^{1-\varepsilon}/\exp{(c_2y(n_{i+1}))}.
\end{equation}
 
\end{thm}

We obtain \eqref{eq1} from \eqref{8} when \eqref{5} holds since in this case
$$
\exp{(c_2y(n))} = o(n).
$$
\section{Preliminary lemmas}

For any non-zero rational number $\alpha$ we may write $\alpha=a/b$ with $a$ and $b$ coprime integers and with $b$ positive. We define $H(\alpha)$, the height of $\alpha$, by
$$
H(\alpha)=\max ( |a|, |b|).
$$
Let $n$ be a positive integer and let $\alpha_1, ..., \alpha_n$ be positive rational numbers with heights at most $A_1,...,A_n$ respectively. Suppose that $A_i \geq 3$ for $i=1,...,n$ and that $\log \alpha_1, ..., \log \alpha_n$ are linearly independent over the rationals where $\log$ denotes the principal value of the logarithm. Let $b_1,...,b_n$ be non-zero integers of absolute value at most $B$ with $B\geq 3$ and put
$$
\Lambda = b_1\log \alpha_1 + ...+ b_n\log \alpha_n.
$$

\begin{lem} \label{lem 1}
There exists an effectively computable positive number $c_0$ such that
$$
\log |\Lambda| > -c_0^n\log A_1 ...\log A_n\log B.
$$
\begin{proof}
This follows from Theorem 2.2 of Nesterenko \cite{N}, which is a special case of the work of Matveev \cite{Matveev}, \cite{Matveev2}.
\end{proof}
\end{lem}

Let $x$ and $y$ be positive real numbers with $y\geq 2$ and let $\Psi (x,y)$ denote the number of positive integers of size at most $x$ all of whose prime factors are of size at most $y$. Let $r$ denote the number of primes of size at most $y$ so that $r=\pi (y).$
\begin{lem} \label{lem 2}
For $2\leq y \leq (\log x)^{1/2}$ we have
$$
\Psi (x,y) = \frac{(\log x)^r}{\prod^{r}_{i=1} (i\log p_i)} (1 +O(y^2(\log x)^{-1}(\log y)^{-1})).
$$

\begin{proof}
This is Theorem 1 of  \cite{Ennola}.
\end{proof}
\end{lem}

We also recall the abc conjecture of Oesterl\'e and Masser \cite{Masser}, \cite{Masser1}, \cite{SY}. Let $x,y$ and $z$ be positive integers. Denote the greatest square-free factor of $xyz$ by $G=G(x,y,z)$ so
$$
G=\prod_{{\substack{p|xyz \\ p , prime}}}p.
$$

\begin{con} \label{con}
(abc conjecture) For each positive real number $\varepsilon$ there is a positive number $c(\varepsilon)$, which depends on $\varepsilon$ only, such that for all pairwise coprime positive integers $x,y$ and $z$ with
$$
x+y=z
$$
we have
$$
z< c(\varepsilon)G^{1+\epsilon}.
$$
\end{con}

For a refinement of the abc conjecture see \cite{RST}.

\section{Proof of Theorem 1}
Let $c_1,c_2,...$ denote effectively computable positive numbers. Following \cite{Tijd}, for $i\geq 3$ we have $n_i\geq 3,$
\begin{equation} \label{9}
n_{i+1} -n_{i} = n_{i}(\frac{n_{i+1}}{n_i} -1)
\end{equation}
and, since $e^z-1 > z$ for $z$ positive,
\begin{equation} \label{10}
\frac{n_{i+1}}{n_i} -1 > \log \frac{n_{i+1}}{n_i}.
\end{equation}
Let $p_1,...,p_r$ be the primes of size at most $y(n_{i+1})$. Notice that $r\geq 2$ since $y(n_{i+1})\geq 3$. 
Then $\frac{n_{i+1}}{n_i} = p_1^{l_1}...p_r^{l_r}$ with $l_1,...,l_r$ integers of absolute value at most $c_1\log n_{i+1}$ and, since $n_{i+1}\leq 2n_i,$
\begin{equation} \label{11}
\max (|l_1|,...,|l_r|) \leq c_2\log n_i.
\end{equation}
Since
$$
\log \frac{n_{i+1}}{n_i} = l_1\log p_1 +\cdots + l_r\log p_r
$$
it follows from \eqref{11} and Lemma \ref{lem 1} that
\begin{equation} \label{12}
\log \frac{n_{i+1}}{n_i}  > (\log n_i)^{-c_3^r\log p_1...\log p_r}.
\end{equation}
By the arithmetic-geometric mean inequality
\begin{equation} \label{13}
\prod^{r}_{i=1} \log p_i \leq (\frac{1}{r}\sum^{r}_{i=1} \log p_i)^r
\end{equation}
and by the prime number theorem
\begin{equation} \label{14}
\sum^{r}_{i=1} \log p_i < c_4r\log r.
\end{equation}
Thus, from \eqref{12}, \eqref{13} and \eqref{14},
\begin{equation} \label{15}
\log \frac{n_{i+1}}{n_i} > (\log n_i)^{-(c_5\log r)^r}.
\end{equation}
Observe that $r\geq 2$ and so
\begin{equation} \label{16}
(c_5\log r)^r < e^{c_6r\log_2 r}.
\end{equation}
Further
$$
3\leq p_r \leq y(n_{i+1})
$$
and so
\begin{equation} \label{17}
r \leq c_7y(n_{i+1})/\log y(n_{i+1}).
\end{equation}
Thus, by \eqref{16} and \eqref{17},
\begin{equation} \label{18}
(c_5\log r)^r < \delta(c_8y(n_{i+1}))
\end{equation}
and \eqref{6} follows from \eqref{9}, \eqref{10}, \eqref{15} and \eqref{18}.

We shall now establish \eqref{7}. Observe that if $n_i$ satisfies \eqref{3} then since $y(t)\geq 3$ for all positive real numbers $t$, $P(2n_i) \leq y(n_i)\leq y(2n_i)$ and so $2n_i=n_j$ for some integer $j$ with $j>n$. In particular $n_{i+1}\leq 2n_i$ hence $n_{i+1}-n_i \leq n_i$ so
\begin{equation} \label{a}
n_{i+1}-n_i < 2n_i.
\end{equation}
Suppose that $X$ is a real number with $X\geq 9$ and that $i$ is a positive integer with $n_{i+1}$ and $n_i$ in the interval $(\sqrt X,X]$. If, in addition,
\begin{equation} \label{b}
y(\sqrt X) > (\log X)^{\frac{1}{4}}
\end{equation}
then, since $\sqrt X < n_i \leq X,$
\begin{equation} \label{d}
y(n_i) > (\log n_i)^{\frac{1}{4}}.
\end{equation}
Since $y$ is non-decreasing
\begin{equation} \label{e}
\pi(y(\sqrt n_i))-1 \leq \pi(y(n_i))
\end{equation}
and by the prime number theorem
$$
\pi(y(n_i)) < c_9\frac{y(n_i)}{\log y(n_i)}.
$$
By \eqref{d},
\begin{equation} \label{f}
\pi(y(n_i)) < c_{10}\frac{y(n_i)}{\log_2 n_i}.
\end{equation}
Thus by \eqref{e} and \eqref{f},
\begin{equation} \label{g}
(\log n_i)^{\pi(y(\sqrt n_i))-1} < e^{c_{10}y(n_i)}.
\end{equation}
We may suppose that $c_1$ exceeds $1+c_{10}$ and in this case, by \eqref{g},
$$
\exp(c_1y(n_i))/(\log n_i)^{\pi(y(\sqrt n_i))-1} \geq \exp(y(n_i)) \geq \exp(3) \geq 2,
$$
and therefore \eqref{7} follows from \eqref{a}.

We shall now show that there is a positive number $c_{11}$ such that if $X$ is a real number with $X> c_{11}$ then there is a positive integer $i$ for which $n_{i+1}$ and $n_{i}$ are in $(\sqrt X,X]$ and satisfy \eqref{7}.
Accordingly let $X$ be a real number with $X\geq 9$ and put
$$
r=\pi(y(\sqrt X)).
$$
Notice that $r\geq 2$ since $y(t)\geq 3$ for all positive real numbers $t$.
By the preceding paragraph we may suppose that
$$
y(\sqrt X) \leq (\log X)^{\frac{1}{4}}.
$$

Let $A(X)$ be the set of integers $n$ with
\begin{equation} \label{19a}
\sqrt X < n \leq X
\end{equation}
for which
\begin{equation} \label{19}
P(n)\leq y(\sqrt X).
\end{equation}
Note that the members of $A(X)$ occur as terms in the sequence $(n_1,n_2,...)$.
The cardinality of $A(X)$ is
$$
\Psi(X,y(\sqrt X)) - \Psi(\sqrt X, y(\sqrt X))
$$
and so for $X > c_{12}$ is, by Lemma \ref{lem 2}, at least
\begin{equation} \label{20}
\frac{(\log X)^r}{2\prod^{r}_{i=1}i\log p_i}.
\end{equation}
Let $j$ be the positive integer for which
$$
\frac{X}{2^j} < \sqrt X \leq \frac{X}{2^{j-1}}
$$
and consider the intervals $(\frac{X}{2^k}, \frac{X}{2^{k-1}}]$ for $k=1,...,j.$
Then $j\leq 1+ \frac{\log X}{2\log 2}$ and so, for $X> c_{13}$,
\begin{equation} \label{21}
j \leq\log X.
\end{equation}
Thus, by \eqref{20} and \eqref{21}, there is an integer $h$ with $1\leq h \leq j$ for which the interval $(\frac{X}{2^h}, \frac{X}{2^{h-1}}]$ contains at least
$$
\frac{(\log X)^{r-1}}{2\prod^{r}_{i=1}i\log p_i}
$$
integers from $A(X)$. Notice that
$$
\prod^{r}_{i=1}i\log p_i  \leq (r \log y(\sqrt X))^r.
$$
Thus, since $y(\sqrt X) \leq (\log X)^{\frac{1}{4}} $ and, since $r \geq 2$, $r-1 \geq \frac{r}{2}$ we see that for $X > c_{14},$ the interval $(\frac{X}{2^h}, \frac{X}{2^{h-1}}]$ contains at least
$$
\frac{(\log X)^{r-1}}{3(r\log y(\sqrt X))^r} + 1
$$
terms from $A(X)$ hence two of them, say $n_{i+1}$ and $n_i$, satisfy
$$
n_{i+1}-n_i < \frac{X}{2^h(\log X)^{r-1}}3(r\log y(\sqrt X))^r.
$$
Since $n_i > \frac{X}{2^h}$ it follows that
$$
n_{i+1}-n_i <  3\frac{n_i}{(\log n_i)^{r-1}}(r\log y(\sqrt X))^r.
$$
By \eqref{19a}, $\sqrt n_i \leq \sqrt X \leq n_i$ hence, since $y$ is non-decreasing, $y(\sqrt n_i) \leq y(\sqrt X) \leq y(n_i)$. Thus
$$
n_{i+1}-n_i < 3\frac{n_i}{(\log n_i)^{r-1}}(r\log y(n_i))^r
$$
and so
\begin{equation} \label{22}
n_{i+1} -n_{i} < 3\frac{n_i}{(\log n_i)^{r'-1}}(s\log y(n_i))^s
\end{equation}
where $r'=\pi(y(\sqrt n_i))$ and $s=\pi(y(n_i))$.
By the prime number theorem there is a positive number $c_{15}$ such that
\begin{equation} \label{23}
3(s\log y(n_i))^s < e^{c_{15}y(n_i)}.
\end{equation}
Estimate \eqref{7} now follows from \eqref{22} and \eqref{23}. On letting $X$ tend to infinity we find infinitely many pairs of integers $n_{i+1}$ and $n_{i}$ which satisfy \eqref{7}.

\section{Proof of Theorem 2}
Let $i\geq 1$ and put
\begin{equation} \label{24}
n_{i+1} -n_{i} = t.
\end{equation}
Let $g$ be the greatest common divisor of $n_{i+1}$ and $n_i$.
Then
$$
\frac{n_{i+1}}{g} - \frac{n_i}{g} = \frac{t}{g}.
$$
Let $\varepsilon > 0$. By the abc conjecture there is a positive number $c(\varepsilon)$ such that
$$
\frac{n_i}{g} < c(\varepsilon)(\frac{t}{g}\prod_{p\leq y(n_{i+1})}p)^{1+\varepsilon}
$$
hence
\begin{equation} \label{25}
(\frac{n_i}{c(\varepsilon)})^{\frac{1}{1+\varepsilon}} < t\prod_{p\leq y(n_{i+1})}p.
\end{equation}
By the prime number theorem, since $y(n_{i+1})\geq 3,$ there exists a positive number $c_2$ such that
\begin{equation} \label{26}
\prod_{p\leq y(n_{i+1})}p < e^{c_2y(n_{i+1})}.
\end{equation}
The result follows from \eqref{24}, \eqref{25} and \eqref{26}.

\section{Acknowledgements}

This research was supported in part by the Canada Research Chairs Program and by grant A3528 from the Natural Sciences and Engineering Research Council of Canada.

\end{document}